\documentclass[a4paper]{ifacconf}
\usepackage{graphicx,amsmath,url}
\usepackage[round]{natbib}
\usepackage{booktabs}
\usepackage{multirow}
\usepackage{amssymb}

\usepackage[T1]{fontenc}
\usepackage[utf8]{inputenc}
\usepackage{ae}
\usepackage{adjustbox}   

{ }
{ }
{ }

\begin{document}

\begin{frontmatter}

\title{Day-ahead Coordination of Virtual Power Plants within
Active Distribution Networks using Deterministic Bi-Level Optimization\thanksref{footnoteinfo}}

\thanks[footnoteinfo]{This work was financially
supported by the Brazilian National Council for Scientific and
Technological Development — CNPq, under project grant 447271/2024-5.}

\author[First]{Laura M. Barajas-Arguello}
\author[Second]{Rafael A. Núñez-Rodríguez}
\author[First]{Daniel Gebbran}
\author[First]{Clodomiro Unsihuay-Vila}

\address[First]{Graduate program in electrical engineering,
Universidade Federal do Paraná, Curitiba, PR, Brazil
(e-mail: laurabarajas@ufpr.br, daniel.gebbran@ufpr.br, clodomiro.vila@ufpr.br).}
\address[Second]{Department of Mechatronic Engineering,
Universidad Santo Tomás, Bucaramanga, SAN, Colombia
(e-mail: rafael.nunez@ustabuca.edu.co).}

\renewcommand{\abstractname}{{\bf Abstract:~}}
\begin{abstract}
This paper proposes a deterministic bilevel optimization
framework for the coordinated operation of Virtual Power
Plants~(VPPs) embedded in an active distribution network.
The Distribution System Operator~(DSO) acts as the
upper-level leader, minimizing a weighted combination of
expenditure, active losses and voltage deviation subject
to nonlinear AC power flow constraints, while each VPP
operates as a lower-level follower that maximizes its
profit under the uniform price signal issued by the DSO.
Unlike most existing formulations, which linearize the
lower-level subproblem to obtain a Mixed-Integer Linear
Program~(MILP), the proposed model retains the full
AC~Optimal Power Flow~(AC-OPF) equations, producing a
bilevel Mixed-Integer Nonlinear Program~(MINLP). The
lower-level problem is replaced by its
Karush-Kuhn-Tucker~(KKT) optimality conditions and the
Strong Duality Theorem, yielding a single-level
Mathematical Program with Equilibrium Constraints~(MPEC).
Complementarity conditions are then linearized via the
Fortuny-Amat big-M transformation. The framework is
validated on the IEEE 33-bus feeder over a 24-hour
horizon, with four distributed resources aggregated into
a single VPP. Compared with individual dispatch against a
regulated time-of-use tariff, aggregation reduces active
losses by 10.5\,\%, the accumulated voltage deviation by
18.3\,\%, and the bus-hours below 0.95\,p.u.\ from 132 to
29. These gains cost 0.31\,\% in social cost and
0.26\,\% in DSO expenditure, while the rent of the
aggregator is preserved.
\end{abstract}

\begin{keyword}
Virtual power plants; Bi-level optimization; AC optimal
power flow; KKT conditions; MPEC; MINLP; Active
distribution network.
\end{keyword}

\end{frontmatter}

\section{Introduction}

Distribution networks have changed more in the last decade
than in the previous fifty. Widespread solar and wind
generation, rising electric vehicle adoption, and battery
storage have turned what was essentially a passive grid
into a system with power injections at many points,
bidirectional flows, and a level of uncertainty that
traditional dispatch was never built to handle
\citep{Gough2023,Zhang2023}. This is exactly the problem
that Virtual Power Plants~(VPPs) are designed for: rather
than managing each distributed resource independently,
they aggregate them as a single entity that can trade with
the market and with the Distribution System
Operator~(DSO) \citep{Zhang2017,Li2022}.

Coordinating a DSO and a VPP is not a minor technical
detail. Both parties have different objectives and do not
share all their information. The DSO wants to minimize
its total operating cost; the VPP wants to maximize its
profit. Treating this as a joint optimization problem,
where a single agent makes all decisions, erases that
asymmetry and biases the solution \citep{Yi2020}.
Stackelberg game theory offers a more realistic
alternative: the DSO acts as the leader and issues a
price signal; the VPP acts as the follower and responds
with its optimal dispatch. The result is a bilevel
problem \citep{Zhang2023,Jadidoleslam2025}.

Representing the network accurately matters.
Most approaches in the literature linearize the power
flow equations to obtain a MILP \citep{Yi2020,Zhang2023,
Vafa2025}. That simplifies computation, but in medium-
voltage distribution networks linear approximations fail
to capture voltage behavior and reactive power losses~---
errors that can exceed 10\,\% in some
cases \citep{Jadidoleslam2025}. This paper keeps the
full nonlinear AC equations at both levels.

Table~\ref{tab:review} summarizes the closest related
work. None of these combines a bilevel structure with an
exact nonlinear AC-OPF lower level. The last column
reports the algebraic modeling language (AML) and the
solver adopted in each study. The survey covers
stochastic models \citep{Gough2023,Jadidoleslam2025},
robust formulations \citep{Vafa2025,Li2022}, distributed
approaches \citep{Li2022}, and unit-commitment-based
frameworks, which confirms that the
gap addressed here is not specific to any single
modeling direction. The closest case, \citet{Yi2020},
uses the same Stackelberg structure with a big-M
reformulation but linearizes the lower level through a DC
approximation; in feeders with high R/X ratios this
discards the voltage and reactive-power behavior that the
AC model retains.

\begin{table*}[t]
\caption{Comparison of related work on VPP optimization.}
\label{tab:review}
\begin{center}
\renewcommand{\arraystretch}{1.15}
\footnotesize
\setlength{\tabcolsep}{2pt}
\begin{tabular}{@{}lcccccccccccccc@{}}
\toprule
\multirow{2}{*}{\textbf{Reference}} &
\multicolumn{2}{c}{\textbf{VPP type}} &
\multicolumn{3}{c}{\textbf{Uncertainty}} &
\multicolumn{3}{c}{\textbf{Control strategy}} &
\multicolumn{3}{c}{\textbf{Optimization model}} &
\multicolumn{2}{c}{\textbf{Solution method}} &
\textbf{AML/Solver} \\
\cmidrule(lr){2-3}\cmidrule(lr){4-6}\cmidrule(lr){7-9}%
\cmidrule(lr){10-12}\cmidrule(lr){13-14}\cmidrule(l){15-15}
& Single & Multi
& Det. & Stoch. & Rob.
& Cent. & Decent. & Dist.
& Bilev. & MILP & MINLP
& KKT & Other & --- \\
\midrule
\citet{Zhang2017}
  &x& & &x& &x& & &x& & & &x& MATLAB / n.r. \\
\citet{Yi2020}
  & &x& & & &x& & &x&x& &x& & GAMS / CPLEX \\
\citet{Wu2021}
  & &x&x& & &x& & & & & & &x& MATLAB / Ipopt \\
\citet{Li2022}
  & &x&x& &x& & &x& & & & &x& MATLAB / Ipopt \\
\citet{Gough2023}
  & &x& &x& &x& & &x&x& & &x& MATLAB / n.r. \\
\citet{Zhang2023}
  & &x& &x& &x& & &x&x& &x& & MATLAB / n.r. \\
\citet{Vafa2025}
  &x& & & &x&x& & &x&x& & &x& GAMS / CPLEX \\
\citet{Jadidoleslam2025}
  & &x& &x& &x& & &x& & &x&x& GAMS / n.r. \\
\citet{Zare2026}
  & &x& &x& &x& & &x&x& & & & GAMS / CPLEX \\
\textbf{This work}
  &x& &x& & &x& & &x& &x&x&x& Pyomo / Ipopt, Bonmin \\
\bottomrule
\end{tabular}
\\[2pt]
\parbox{\textwidth}{\scriptsize \emph{n.r.}: not reported.
MATLAB is listed as the modeling environment when the
study does not specify an AML.}
\end{center}
\end{table*}

The paper makes three contributions. First, the full
nonlinear AC-OPF equations are kept in the lower level~---
no DC approximation, no LinDistFlow. Second, the
Stackelberg structure preserves the economic independence
of the DSO and the VPP, which joint-optimization models
cannot guarantee. Third, the bilevel problem is converted
into a tractable MINLP through KKT conditions and the
Fortuny-Amat big-M transformation. Results compare
Case~1, in which the resources are dispatched
individually against the regulated time-of-use tariff,
and Case~2, in which the same resources are aggregated
into a VPP governed by the bilevel model, on the IEEE
33-bus feeder over a 24-hour horizon.

\section{Bilevel Problem Formulation}

The interaction between the DSO and a set of
aggregators $v\in\mathcal{V}$ follows a Stackelberg
structure. The DSO (leader) issues a uniform price
signal $\lambda_{t}$ and each VPP (follower) responds
with the optimal dispatch of its own resources
\citep{Zhang2023,Yi2020}. The upper-level decision
vector is $\mathbf{x}^{U}=\{P_{k,t}^G, Q_{k,t}^G,
V_{n,t}, \delta_{n,t}, \lambda_{t}\}$ and the
lower-level vector is $\mathbf{x}^{L}=\{P_{t}^{VPP},
P_{r,t}\}$. The formulation below is written for a
single aggregator, $|\mathcal{V}|=1$, which is the
configuration validated in Section~4; extending it to
$|\mathcal{V}|>1$ replicates the KKT block per
aggregator without altering the reformulation.
All sets, parameters, and variables are defined in
Appendix~A. Fig.~\ref{fig:estructura} illustrates the
model structure.

\begin{figure}[ht]
\begin{center}
\includegraphics[width=8.4cm]{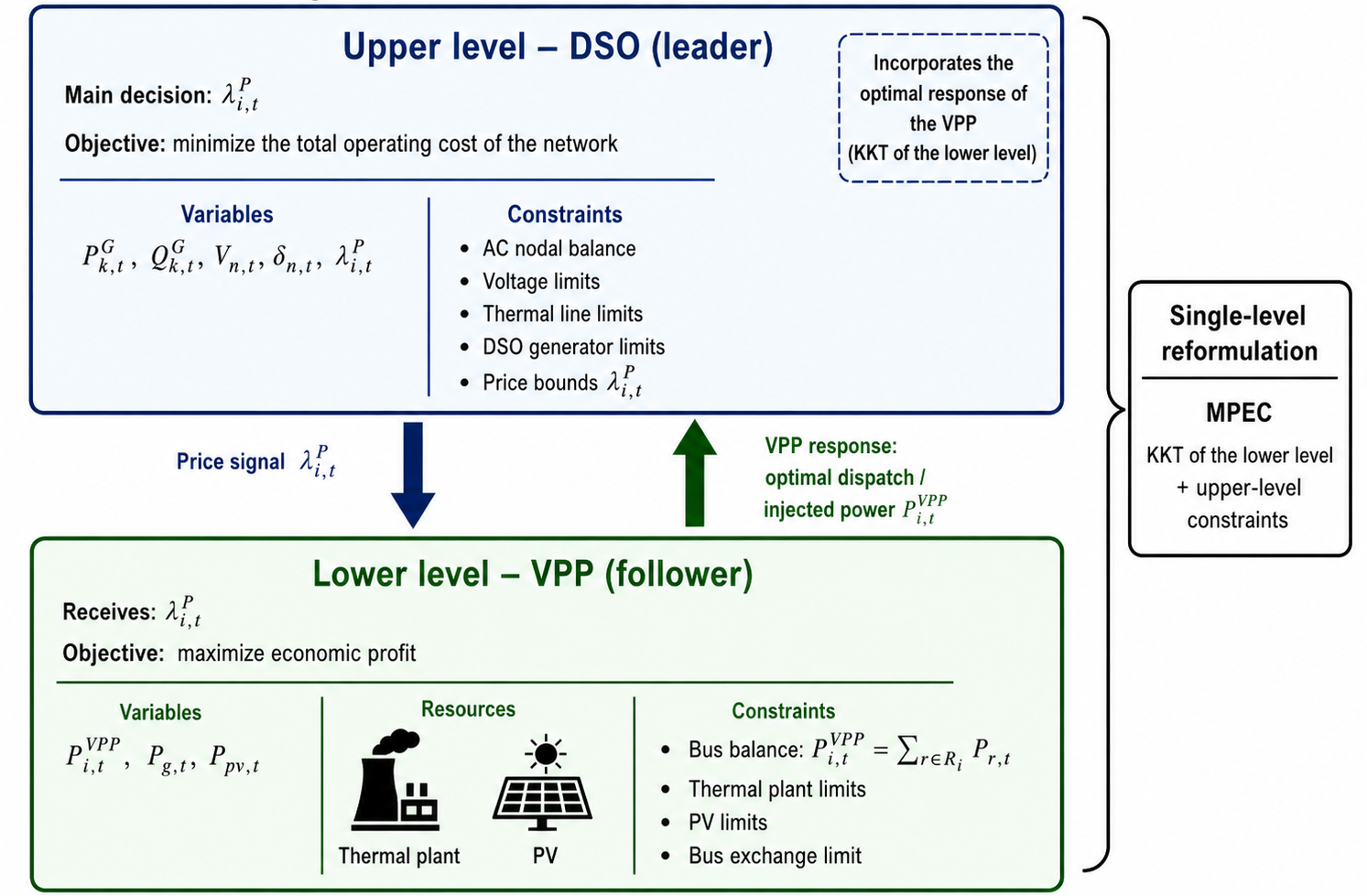}
\caption{Structure of the DSO--VPP bilevel model. The
upper level (DSO) issues $\lambda_{t}$; the lower level
(VPP) responds with the optimal dispatch of its resources.}
\label{fig:estructura}
\end{center}
\end{figure}

\subsection{Lower level -- VPP (follower)}

Given $\lambda_{t}$, the VPP maximizes its net
economic benefit over the 24-hour horizon $T$:
\begin{equation}
\max_{\mathbf{x}^{L}}\; B^{VPP} =
  \sum_{t\in T}\!
  \bigl[\lambda_{t} P_{t}^{VPP}
        - C_g(P_{g,t}) - C_{pv}(P_{pv,t})\bigr]
\label{eq:obj_lower}
\end{equation}
with operating costs:
\begin{align}
C_g(P_{g,t})     &= a_g P_{g,t}^{2}
                    + b_g P_{g,t} + c_g,\; a_g>0
                    \label{eq:cost_th}\\
C_{pv}(P_{pv,t}) &= \pi_{pv} P_{pv,t},\;
                    \pi_{pv}\geq 0
                    \label{eq:cost_pv}
\end{align}

The quadratic cost for the thermal unit is standard in
VPP economic dispatch \citep{Jadidoleslam2025}. The
photovoltaic (PV) unit carries only a linear O\&M cost
since solar fuel is free \citep{Zhang2023}.
The lower-level constraints are:
\begin{align}
P_{t}^{VPP} &= \textstyle\sum_{r\in R} P_{r,t},
  &\forall t\!\in\! T \label{eq:ll_balance}\\
P_g^{\min} &\leq P_{g,t} \leq P_g^{\max},
  &\forall t\!\in\! T \label{eq:ll_th}\\
0 &\leq P_{pv,t} \leq \bar{P}_{pv,t},
  &\forall t\!\in\! T \label{eq:ll_pv}\\
0 &\leq P_{t}^{VPP} \leq P^{VPP,\max},
  &\forall t\!\in\! T \label{eq:ll_vpp}\\
\textstyle\sum_{t\in T} &P_{g,t}\,\Delta t \leq E_g^{\max},
  &\forall g \label{eq:ll_energy}
\end{align}
where $\bar{P}_{pv,t}$ is the available solar power at
period $t$ (irradiance profile, external parameter), and
$R$ collects the thermal and photovoltaic units of the
aggregator, so the VPP injection is the sum of their
contributions.
Constraint~(\ref{eq:ll_energy}) caps the daily energy
produced by each thermal unit at $E_g^{\max}$, with
$\Delta t = 1$\,h. It represents the fuel budget available
to the aggregator over the horizon. Without it, the
thermal units operate near their rated power throughout
the day whenever the price signal exceeds their marginal
cost, which is not representative of a distributed
generation unit.

The follower's feasible set contains only its internal
resource limits and the exchange limit at its
connection point. Voltage magnitudes, line flows and
the AC balance are enforced exclusively at the upper
level, where $P^{VPP}_{t}$ enters
Eq.~(\ref{eq:pf_p}) as a nodal injection. This
reflects the information asymmetry that the
Stackelberg structure is meant to represent: the
aggregator does not hold the network data required to
evaluate those constraints. Network feasibility is
nevertheless guaranteed by construction, because the
single-level reformulation enforces the AC constraints
and the follower optimality conditions simultaneously,
so no price-response pair that violates the network
can be part of the equilibrium.

\subsection{Upper level -- DSO (leader)}

The DSO minimizes a weighted sum of three normalized
objectives: its total expenditure, the active energy
losses of the feeder, and the accumulated voltage
deviation \citep{Yi2020,Zhang2023}:
\begin{equation}
\min_{\mathbf{x}^{U}}\; F =
  W_{C}\frac{C^{DSO}}{C^{ref}}
+ W_{L}\frac{P^{loss}}{P^{loss,ref}}
+ W_{V}\frac{F^{V}}{F^{V,ref}}
\label{eq:obj_upper}
\end{equation}
with
\begin{equation}
C^{DSO} = \sum_{t\in T}\Bigl[\sum_{k\in G}
  \bigl(a_k (P_{k,t}^G)^{2}+ b_k P_{k,t}^G + c_k\bigr)
  + \lambda_{t}P_{t}^{VPP}\Bigr]
\label{eq:c_dso}
\end{equation}
where where $P^{loss}=\sum_{t}\sum_{(n,m)\in L}G_{nm}
(V_{n,t}^{2}+V_{m,t}^{2}\allowbreak
{}-2V_{n,t}V_{m,t}\cos\theta_{nm,t})$
is the total active loss of the feeder and 
$F^{V}=\sum_{t}\sum_{n\in N}(V_{n,t}-1)^{2}$
is the accumulated squared deviation from nominal
voltage.

The three terms are normalized against the
corresponding Case~1 values $C^{ref}$, $P^{loss,ref}$
and $F^{V,ref}$, so $F=1$ reproduces the
non-aggregated operation and $F<1$ denotes a joint
improvement. The weights are $W_C=0.6$ and
$W_L=W_V=0.2$, which keeps cost as the dominant
criterion while allowing the price signal to be shaped
by network conditions.

The constraints include the nonlinear AC power flow
equations for every bus $n\!\in\!N$ and period $t\!\in\!T$:
\begin{align}
&\sum_{k\in G_n}\!P_{k,t}^G
  +\!\sum_{r\in R_n}\!P_{r,t} - P_{n,t}^D \notag\\
&= \sum_{m\in N}\!V_{n,t}V_{m,t}
  (G_{nm}\cos\theta_{nm,t}
  +B_{nm}\sin\theta_{nm,t}) \label{eq:pf_p}\\[4pt]
&\sum_{k\in G_n}\!Q_{k,t}^G - Q_{n,t}^D \notag\\
&= \sum_{m\in N}\!V_{n,t}V_{m,t}
  (G_{nm}\sin\theta_{nm,t}
  -B_{nm}\cos\theta_{nm,t}) \label{eq:pf_q}
\end{align}
where $\theta_{nm,t}=\delta_{n,t}-\delta_{m,t}$.
Active and reactive line flows are:
\begin{multline}
P_{nm,t} = V_{n,t}^2 G_{nm} - V_{n,t}V_{m,t}\\
  \cdot(G_{nm}\cos\theta_{nm,t}
  + B_{nm}\sin\theta_{nm,t}) \label{eq:pnm}
\end{multline}
\begin{multline}
Q_{nm,t} = -V_{n,t}^2 B_{nm} - V_{n,t}V_{m,t}\\
  \cdot(G_{nm}\sin\theta_{nm,t}
  - B_{nm}\cos\theta_{nm,t}) \label{eq:qnm}
\end{multline}
Operational limits:
\begin{align}
{-}P_{nm}^{\max} &\leq P_{nm,t} \leq P_{nm}^{\max},
  &\forall (n,m)\!\in\!L,\;t\!\in\!T \label{eq:ul_line}\\
P_k^{\min}&\leq P_{k,t}^G \leq P_k^{\max},
  &\forall k\!\in\!G,\;t\!\in\!T \label{eq:ul_pg}\\
Q_k^{\min}&\leq Q_{k,t}^G \leq Q_k^{\max},
  &\forall k\!\in\!G,\;t\!\in\!T \label{eq:ul_qg}\\
V_n^{\min}&\leq V_{n,t} \leq V_n^{\max},
  &\forall n\!\in\!N,\;t\!\in\!T \label{eq:ul_v}\\
{-}2\pi &\leq \delta_{n,t} \leq 2\pi,
  &\forall n\!\in\!N,\;t\!\in\!T \label{eq:ul_delta}\\
\kappa^{P,\min}&\leq \lambda_{t}
                \leq \kappa^{P,\max},
  &\forall t\!\in\!T \label{eq:ul_price}\\
B^{VPP} &\geq B^{ref} \label{eq:ul_ir}
\end{align}

Constraint~(\ref{eq:ul_ir}) is the individual
rationality condition: the aggregator will not accept
a price signal that leaves it worse off than the
outside option $B^{ref}$, taken as the rent it obtains
in Case~1 under the regulated tariff. Written in dual
variables through strong duality,
\begin{multline}
B^{VPP}=\sum_{r,t}a_{r}P_{r,t}^{2}
 +\sum_{r,t}\overline{\nu}_{r,t}\overline{P}_{r,t}\\
 +\sum_{t}\overline{\nu}^{V}_{t}P^{VPP,\max}
 +\sum_{r}\theta_{r}E_{r}^{\max}
\label{eq:benef_dual}
\end{multline}

\section{Solution Procedure}

Bilevel problems are generally intractable when solved
directly, because the lower level appears as an implicit
constraint of the upper level. The approach here converts
the bilevel into a single-level problem in three steps,
following the standard methodology in the VPP
literature \citep{Yi2020,Zhang2023,Gough2023}: the KKT
conditions of the lower level, the strong duality
substitution that yields the MPEC, and the big-M
linearization of the complementarity conditions.

\subsection{KKT conditions of the lower level}

The lower level is convex: objective~(\ref{eq:obj_lower})
is strictly concave in $P_{g,t}$ and linear in all other
variables; all constraints are linear. The KKT conditions
are therefore necessary and sufficient for
optimality \citep{Yi2020,Zhang2023}.

Associating $\mu_{t}$ (unrestricted in sign)
with~(\ref{eq:ll_balance}), and nonnegative multipliers
$\underline{\nu}_{g,t},\overline{\nu}_{g,t}$
with~(\ref{eq:ll_th}),
$\underline{\nu}_{pv,t},\overline{\nu}_{pv,t}$
with~(\ref{eq:ll_pv}),
$\underline{\nu}_{t}^{V},\overline{\nu}_{t}^{V}$
with~(\ref{eq:ll_vpp}) and $\theta_g\geq 0$
with~(\ref{eq:ll_energy}), the stationarity conditions are:
\begin{align}
2a_g P_{g,t}+b_g-\mu_{t}
  -\underline{\nu}_{g,t}+\overline{\nu}_{g,t}+\theta_g
  &= 0 \label{eq:kkt_g}\\
\pi_{pv}-\mu_{t}
  -\underline{\nu}_{pv,t}+\overline{\nu}_{pv,t}
  &= 0 \label{eq:kkt_pv}\\
-\lambda_{t}+\mu_{t}
  -\underline{\nu}_{t}^{V}+\overline{\nu}_{t}^{V}
  &= 0 \label{eq:kkt_vpp}
\end{align}

The complementary slackness conditions are:
\begin{equation}
\begin{adjustbox}{max width=0.92\columnwidth}
$\displaystyle
\begin{aligned}
\underline{\nu}_{g,t}(P_{g,t}-P_g^{\min}) &= 0,\;
\overline{\nu}_{g,t}(P_g^{\max}-P_{g,t}) = 0\\
\underline{\nu}_{pv,t}\,P_{pv,t} &= 0,\;
\overline{\nu}_{pv,t}(\bar{P}_{pv,t}-P_{pv,t}) = 0\\
\underline{\nu}_{t}^{V}\,P_{t}^{VPP} &= 0,\;
\overline{\nu}_{t}^{V}(P^{VPP,\max}-P_{t}^{VPP}) = 0\\
\theta_g\Bigl(E_g^{\max}-\sum_{t\in T}
  P_{g,t}\,\Delta t\Bigr) &= 0
\end{aligned}
$
\end{adjustbox}
\label{eq:cs}
\end{equation}

\subsection{Strong duality and MPEC reformulation}

The term $\lambda_{t}\cdot P_{t}^{VPP}$ is bilinear:
$\lambda_{t}$ is a leader decision variable and
$P_{t}^{VPP}$ is a follower variable. Left untreated,
the single-level problem would remain nonconvex and hard
to solve.

The Strong Duality Theorem resolves this. Since the lower
level is convex, the optimal primal and dual objective
values are equal, $f^*=g^*$ \citep{Yi2020}:
\begin{equation}
f(\mathbf{x}^{L*}) = g(\boldsymbol{\mu}^*,
\underline{\boldsymbol{\nu}}^*,
\overline{\boldsymbol{\nu}}^*)
\label{eq:strong_dual}
\end{equation}

This equality lets the bilinear term be replaced by an
expression in follower variables and Lagrange
multipliers \citep{Yi2020}:
\begin{multline}
\sum_{t\in T}\lambda_{t}P^{VPP}_{t}=
  \sum_{r,t}\bigl[2a_{r}P_{r,t}^{2}+b_{r}P_{r,t}\bigr]
  -\sum_{r,t}\underline{\nu}_{r,t}P_{r}^{\min}\\
  +\sum_{r,t}\overline{\nu}_{r,t}\overline{P}_{r,t}
  +\sum_{t}\overline{\nu}^{V}_{t}P^{VPP,\max}
  +\sum_{r}\theta_{r}E_{r}^{\max}
\label{eq:dual_sub}
\end{multline}
where $\overline{P}_{r,t}$ denotes the upper bound of
resource $r$ at period $t$, equal to $P_g^{\max}$ for
the thermal units and to $\bar{P}_{pv,t}$ for the
photovoltaic units. The multiplier $\mu_{t}$ cancels
when the stationarity conditions are substituted,
which is precisely what makes the identity useful: the
bilinear coupling between the two levels disappears
and only follower variables remain.

Substituting~(\ref{eq:dual_sub}) into the expenditure
term~(\ref{eq:c_dso}) of~(\ref{eq:obj_upper}) converts
the bilevel into a single-level MPEC:
\begin{multline}
C^{DSO}=\sum_{t\in T}\sum_{k\in G}
  \bigl(a_k(P_{k,t}^G)^{2}+b_k P_{k,t}^G+c_k\bigr)\\
  +\sum_{r,t}\bigl[2a_{r}P_{r,t}^{2}+b_{r}P_{r,t}\bigr]
  -\sum_{r,t}\underline{\nu}_{r,t}P_{r}^{\min}\\
  +\sum_{r,t}\overline{\nu}_{r,t}\overline{P}_{r,t}
  +\sum_{t}\overline{\nu}^{V}_{t}P^{VPP,\max}
  +\sum_{r}\theta_{r}E_{r}^{\max}
\label{eq:mpec_obj}
\end{multline}

\subsection{Big-M linearization (MINLP)}

The complementarity conditions
are discontinuous and cannot be handled directly by
standard solvers. The Fortuny-Amat
transformation \citep{FortunyAmat1981,Yi2020}
introduces binary variables
$u\in\{0,1\}$: for each pair $(\nu,\,g(\mathbf{x}))$,
the constraints $0\leq\nu\leq M\cdot u$ and
$0\leq g(\mathbf{x})\leq M(1-u)$ force at least one to
be zero.
The big-M bounds are derived from problem
parameters \citep{Yi2020}. Those associated with primal
slacks take the width of the corresponding feasible
interval, $M_g=P_g^{\max}-P_g^{\min}$,
$M_{pv}=\max_{t}\bar{P}_{pv,t}$, $M_V=P^{VPP,\max}$ and
$M_E=E_g^{\max}$. Those associated with multipliers take
the largest marginal value the dual can reach,
$M_{\nu g}=2a_g P_g^{\max}+b_g+\mu^{\max}$,
$M_{\nu pv}=\pi_{pv}+\mu^{\max}$,
$M_{\nu V}=\kappa^{P,\max}+\mu^{\max}$ and
$M_{\theta}=2a_g P_g^{\max}+b_g$, with
$\mu^{\max}=\max\{\kappa^{P,\max},\,
2a_g P_g^{\max}+b_g,\,\pi_{pv}\}$.

\section{Case Studies and Results}

\begin{figure}[ht]
\begin{center}
\includegraphics[width=8.5cm]{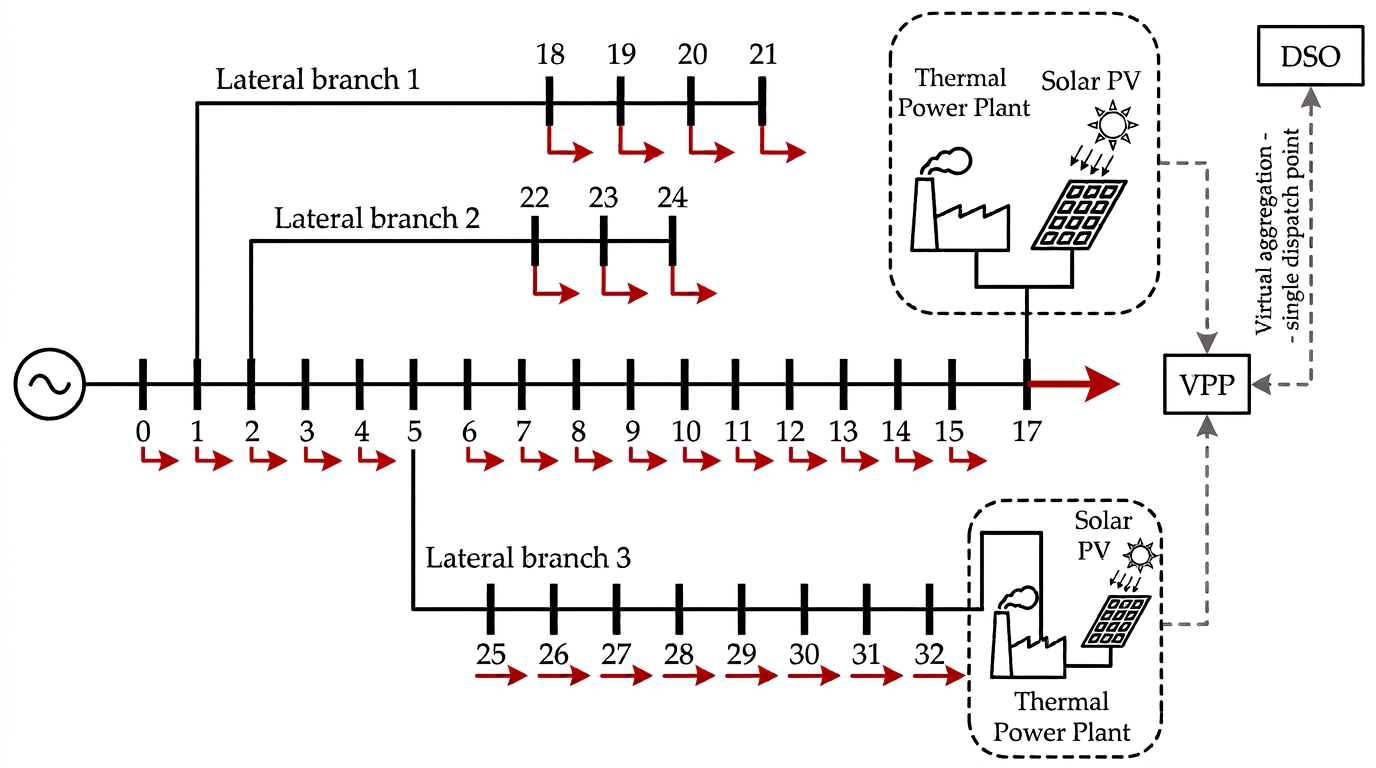}
\caption{IEEE 33-bus radial distribution feeder with
the VPP. The four resources at buses~17 and~32 are
aggregated into a single dispatch point that responds
to the price signal $\lambda_{t}$. In Case~1 the same
four resources occupy the same buses but are
dispatched individually against the regulated
time-of-use tariff, without aggregation. Red arrows
denote bus loads.}
\label{fig:red33b}
\end{center}
\end{figure}

The test system is the IEEE 33-bus radial distribution
feeder \citep{BaranWu1989} (Fig.~\ref{fig:red33b}),
with 32 lines and a base
power of 100\,MVA. Bus~0 is the slack bus; the remaining
32 buses are load buses. Peak demand is 3.715\,MW and
2.300\,MVAr. Both active demand $P_{n,t}^D$ and reactive
demand $Q_{n,t}^D$ follow 24-hour profiles ranging from
0.50 to 1.00 of the peak value, and are enforced through
the AC balance in
Eqs.~(\ref{eq:pf_p})--(\ref{eq:pf_q}). Voltage magnitudes
are bounded to 0.9--1.1\,p.u. and each line has an
apparent power limit of 10\,MVA. The substation generator
at bus~0 has cost coefficients $a_k=25$\,R\$/MWh$^2$,
$b_k=297$\,R\$/MWh, $c_k=0$, with limits of 0 to 10\,MW
and $\pm10$\,MVAr.

The VPP aggregates four resources placed at buses~17
and~32, identified as the electrically weakest nodes of
the feeder through the accumulated series resistance
along the radial path from the substation. Each of the
two buses hosts one thermal unit
($P_g^{\min}=0$, $P_g^{\max}=0.4$\,MW, $a_g=24$\,R\$/MWh$^2$,
$b_g=470$\,R\$/MWh, $c_g=0$, $E_g^{\max}=3.84$\,MWh) and one PV unit
($\bar{P}_{pv}=0.6$\,MW at peak irradiance,
$\pi_{pv}=5$\,R\$/MWh, covering O\&M only), giving
0.8\,MW of dispatchable capacity and 1.2\,MW of solar
capacity. The irradiance profile is a synthetic
half-sinusoidal curve, nonzero from hour~7 to hour~17 and
peaking at noon. The daily energy budget corresponds to a 40\% capacity
factor for each thermal unit over the 24-hour horizon.
The price signal is bounded by
$\kappa^{P,\min}=57.31$ and
$\kappa^{P,\max}=785.27$\,R\$/MWh, corresponding to the
structural PLD floor and ceiling defined by ANEEL for
2026 \citep{ANEEL2025}.

In Case~1 each resource self-dispatches against a
regulated time-of-use tariff whose periods follow the
ANEEL white tariff structure, with a three-hour peak
window from hour~18 to hour~20 and intermediate hours
immediately before and after it \citep{ANEEL2021}.
The tariff levels are scaled from the base-case
average marginal cost of 455.30\,R\$/MWh so that the
daily average matches it, giving 407.37\,R\$/MWh
off-peak, 551.15\,R\$/MWh intermediate and
694.93\,R\$/MWh at peak. The resulting injections
enter the AC-OPF as fixed parameters and the problem
is solved with Ipopt. Case~2 is the single-level MINLP
obtained from the bilevel reformulation, with 2\,210
variables (242 binary) and 2\,983 constraints, solved
with Bonmin.

\subsection{DSO dispatch and cost}

In Case~1 the substation supplies most of the demand
while the four resources self-dispatch against the
regulated time-of-use tariff. Over the 24-hour horizon
the substation generation cost totals
\textbf{22\,399.01\,R\$} and the payment to the
resources \textbf{6\,307.43\,R\$}, of which
1\,963.97\,R\$ corresponds to production cost, giving
a total DSO expenditure of \textbf{28\,706.44\,R\$}.

In Case~2 the aggregation shifts the balance between
these two components. Substation generation cost falls
to \textbf{20\,731.79\,R\$}, a reduction of 7.44\,\%,
while resource production cost rises to
3\,707.21\,R\$ and the payment to the VPP rises to
\textbf{8\,050.66\,R\$}, so DSO expenditure increases
marginally to \textbf{28\,782.45\,R\$}, only
0.26\,\% above Case~1. Expenditure alone is an
incomplete metric, because the payment is a transfer
rather than a consumption of resources. The social
cost, defined as substation generation cost plus
resource production cost and excluding transfers,
amounts to \textbf{24\,362.98\,R\$} in Case~1 and
\textbf{24\,439.00\,R\$} in Case~2, an increase of
0.31\,\%. Both conventions lead to the same conclusion:
preserving the economic autonomy of the aggregator
costs the system less than half a percent, and this is
the price paid for the technical benefits reported
below.

\subsection{VPP dispatch and price signal}

\begin{figure}[t]
\begin{center}
\includegraphics[width=8.5cm]{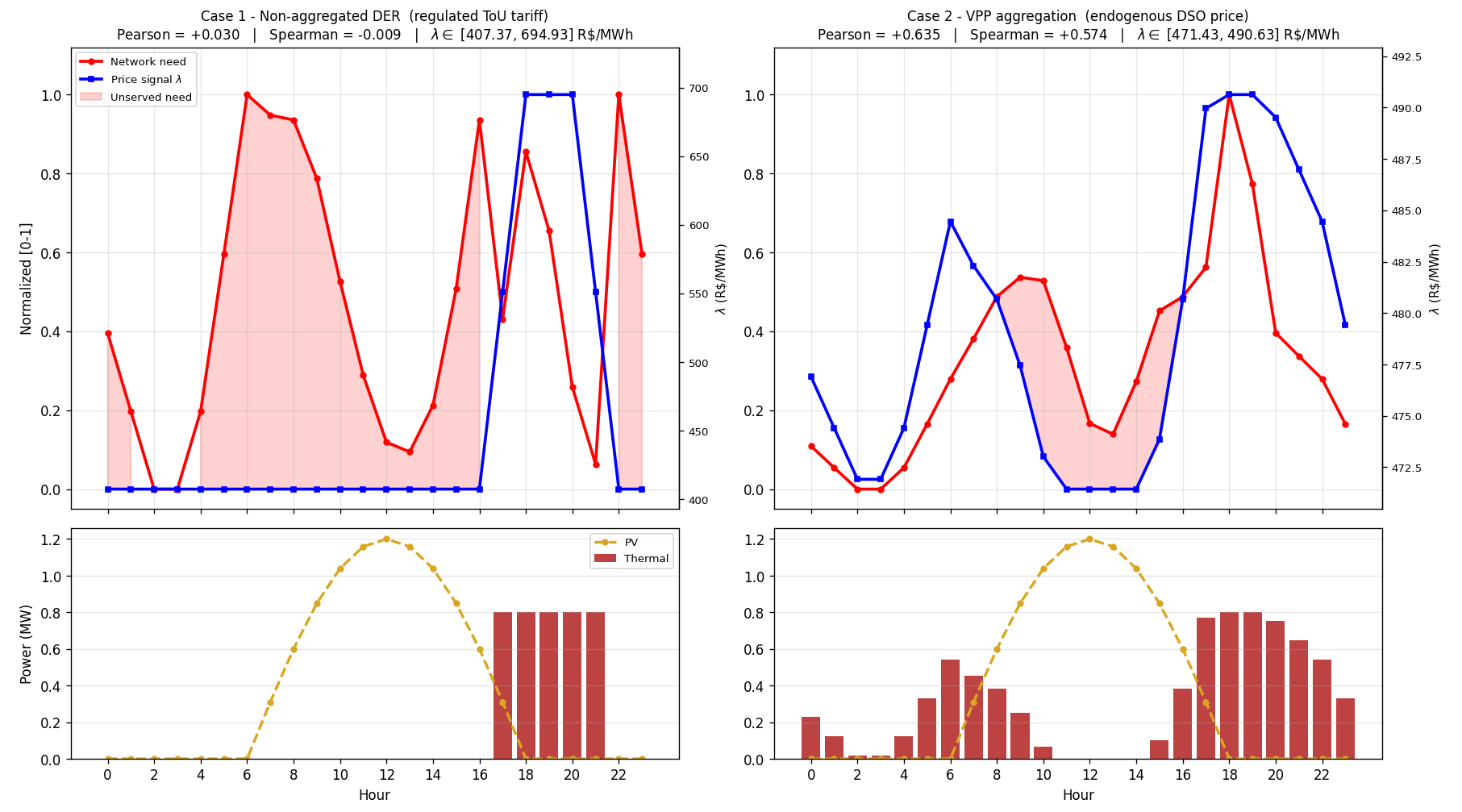}
\caption{Price signal against the network need
indicator. (a)~Case~1, regulated time-of-use tariff.
(b)~Case~2, price signal issued by the DSO. The
left-hand axis is min-max normalized and the
right-hand axis gives the actual price, so a
normalized value of zero corresponds to the off-peak
tariff of 407.37\,R\$/MWh in Case~1 and not to a zero
price. Shaded areas denote network need that the price
signal does not reward. Case~2 prices follow from the
minimum-dispersion dual selection.}
\label{fig:precio}
\end{center}
\end{figure}

The daily energy budget binds for both thermal units,
which reach a 40\,\% capacity factor and deliver
7.680\,MWh over the horizon, against 4.000\,MWh in
Case~1, where they run only during the five
tariff-peak hours. Photovoltaic output is identical in
both cases at 9.115\,MWh with zero curtailment. The
aggregate VPP injection peaks at 1.200\,MW at noon.
The price signal issued by the DSO averages
479.11\,R\$/MWh and stays within a narrow band of
471.43 to 490.63\,R\$/MWh
(Fig.~\ref{fig:precio}). Its correlation with the
network need indicator reaches $\rho=0.635$, against
$\rho=0.030$ for the tariff signal of Case~1. The
aggregator profit is 4\,343.45\,R\$ in both cases, so
the technical gains reported below are obtained
without eroding the economic position of the
aggregator.

The signal is lowest during the solar hours and
highest at the evening peak, when photovoltaic output
has vanished and demand is at its maximum. It tracks
the marginal value of a local injection rather than
the availability of solar energy, which is why
abundant irradiance depresses rather than raises it.
The fuel shadow price $\theta_{r}=1.43$\,R\$/MWh
reallocates the limited thermal energy towards those
same evening hours. The signal is not a substitute for
the wholesale settlement price: the PLD floor and
ceiling bound its admissible range, while
$\lambda_{t}$ differentiates within that range
according to the state of the distribution network.
The two are therefore complementary, and the
correlation reported above supports this reading.

\subsection{Voltage profiles}

\begin{figure}[ht]
\begin{center}
\includegraphics[width=8.5cm]{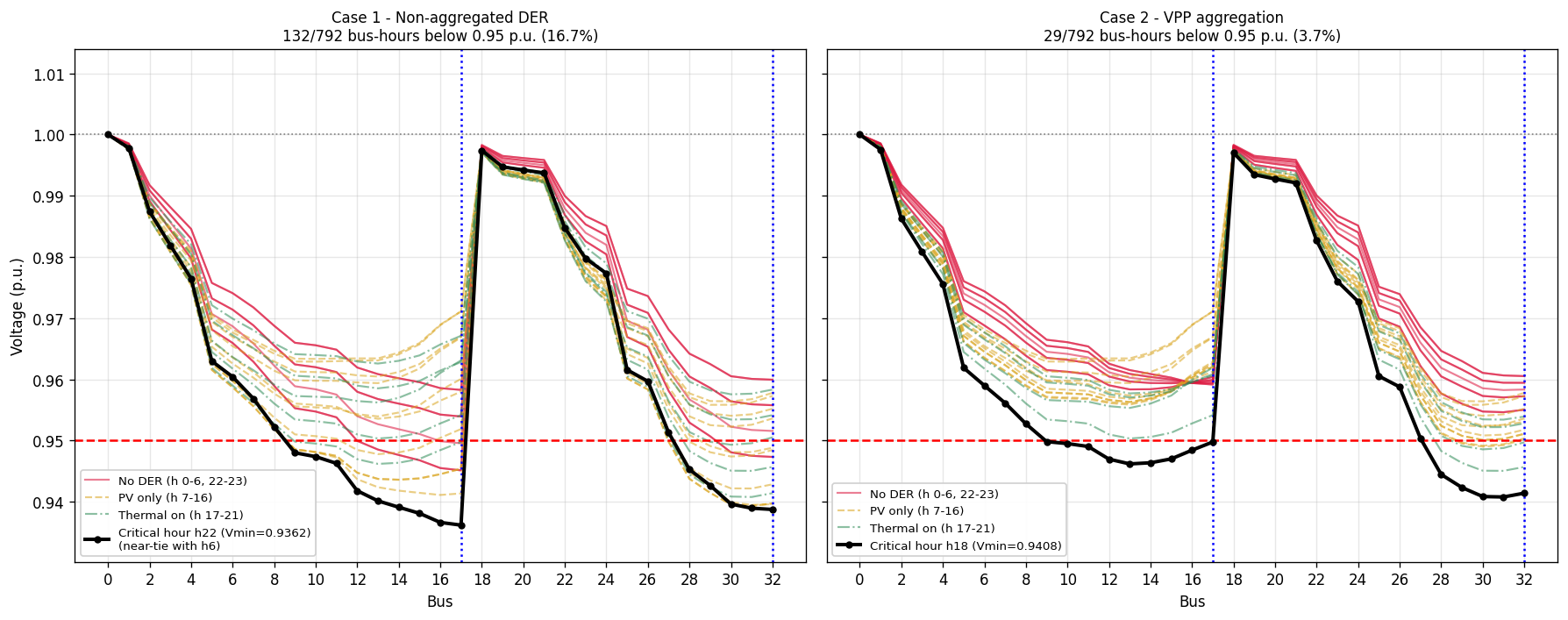}
\caption{Voltage profile along the feeder grouped by
operating regime. (a)~Case~1, non-aggregated
resources. (b)~Case~2, VPP aggregation. Curves are
grouped into hours without resource output, hours with
photovoltaic output only, and hours with thermal
output. The black curve is the critical hour of each
case. Dotted vertical lines mark the resource buses~17
and~32; the dashed red line is the 0.95\,p.u.\
reference.}
\label{fig:tension}
\end{center}
\end{figure}

The 0.95\,p.u.\ threshold is the lower bound of the
adequate voltage range defined by ANEEL for
distribution systems, and is used here as a quality
indicator rather than as a model constraint, since the
optimization enforces the wider $[0.9,\,1.1]$\,p.u.\
range \citep{ANEELPRODIST8}. In Case~1 the minimum
voltage is 0.9362\,p.u.\ at bus~17 during hour~22, and
132 of the 792 bus-hours fall below 0.95\,p.u. In
Case~2 the minimum rises to 0.9408\,p.u.\ at bus~31
during hour~18, and the count drops to 29 bus-hours, a
reduction of 78.0\,\%. The accumulated voltage
deviation $F^{V}$ falls from 1.0448 to
0.8533\,p.u.$^2$, a reduction of 18.3\,\%. Both cases
respect the $[0.9,\,1.1]$\,p.u.\ limits at all times
(Fig.~\ref{fig:tension}).

\subsection{Power flows and losses}

\begin{figure}[ht]
\begin{center}
\includegraphics[width=8.5cm]{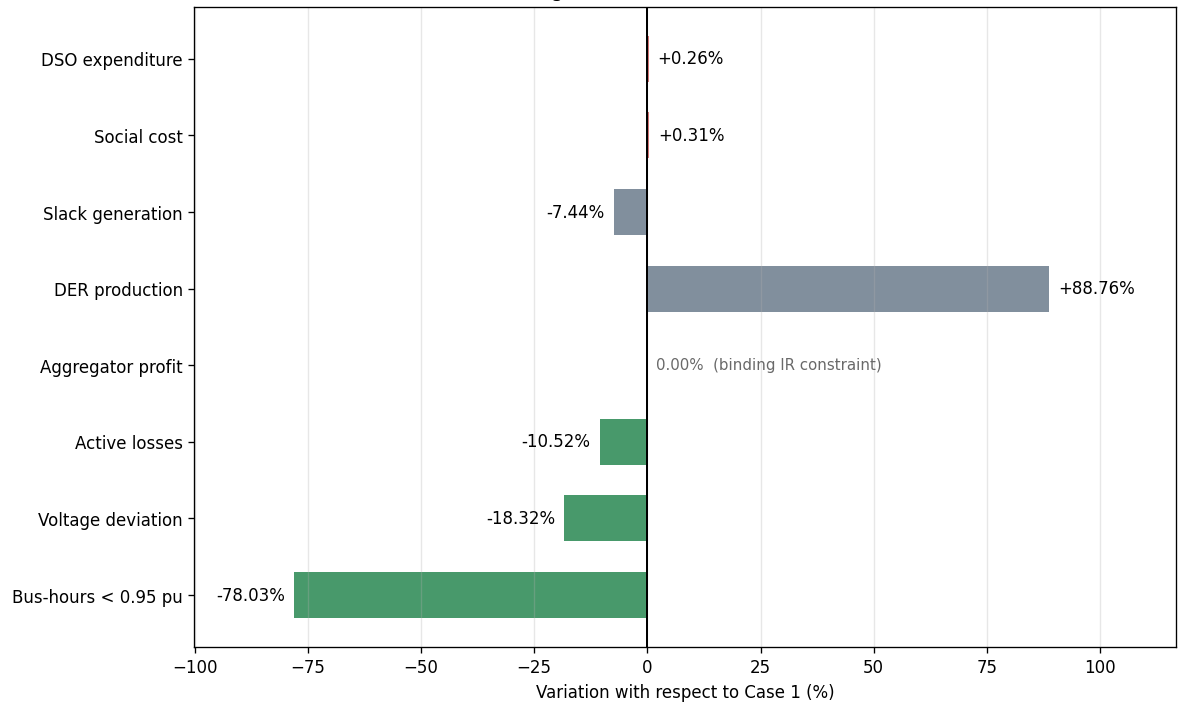}
\caption{Case~2 against Case~1, variation of each
indicator. The grey bars are the cost mechanism rather
than an outcome: local generation displaces
centralised generation and the two nearly cancel. The
aggregator profit is unchanged because the individual
rationality constraint~(\ref{eq:ul_ir}) is binding.}
\label{fig:indicadores}
\end{center}
\end{figure}

Total active losses fall from 2.2293\,MWh in Case~1 to
1.9949\,MWh in Case~2, a reduction of 10.5\,\%. The
gain comes from relocating generation closer to the
weak end of the feeder during the hours of highest
loading, which the tariff signal of Case~1 does not
reward. Fig.~\ref{fig:indicadores} summarizes the
variation of every indicator.

\subsection{KKT verification}

The strong duality identity is satisfied with a
residual of $1.8\times10^{-12}$\,R\$ after the dual
selection stage, and the maximum complementarity
violation across all periods and resources is
$1.2\times10^{-5}$, which validates the big-M
implementation.

\section{Conclusion}

A deterministic bilevel framework for DSO--VPP
coordination was presented, retaining the full
nonlinear AC-OPF in the lower level and reducing it to
a single-level MINLP through KKT conditions, strong
duality and big-M linearization. On the IEEE 33-bus
feeder over 24 hours, aggregation lowers active losses
by 10.5\,\%, the accumulated voltage deviation by
18.3\,\%, and the bus-hours below 0.95\,p.u.\ from 132
to 29, at a cost of 0.31\,\% in social cost and
0.26\,\% in DSO expenditure. The individual rationality
constraint holds the aggregator at its outside option,
so that cost falls on the system and not on a transfer.

Several extensions follow from this work. The
immediate one is the coordination of multiple
aggregators, $|\mathcal{V}|>1$, where the uniform
price signal must be tested against locational
alternatives. Energy storage under coupled
charge/discharge constraints and a smart transformer
as a controllable network asset would both enlarge the
feasible set of the leader. A comparison against a DC
approximation of the lower level would quantify the
accuracy gains of the AC formulation in voltage and
reactive-power terms \citep{Jadidoleslam2025}, and the
framework also admits stochastic uncertainty in
photovoltaic generation and CVaR risk measures.

\section{Agradecimentos}
Acknowledgment of support for Agreement 148/2026 PDI - Agreement for Research,
Development, and Innovation (R, D\,\&I Agreement) between the Araucária Foundation
and the Federal University of Paraná.

\bibliography{ifacconf}

\end{document}